\documentclass[12pt]{amsart}
\input{xy}
\xyoption{all}
\usepackage{amssymb}

\newtheorem{theorem}{Theorem}[section]
\newtheorem{lemma}[theorem]{Lemma}

\newtheorem{proposition}[theorem]{Proposition}
\newtheorem{corollary}[theorem]{Corollary}

\newtheorem{remark}[theorem]{Remark}

\numberwithin{equation}{section}

\newcommand{\ra}{\rightarrow}
\newcommand{\la}{\leftarrow}

\def\umono{\ar@{_{(}->}[u]}
\def\uumono{\ar@{_{(}->}[uu]}

\def\lmono{\ar@{_{(}->}[l]}
\def\llmono{\ar@{_{(}->}[ll]}

\newcommand{\map}{{\rm map}\,}

\newcommand{\Z}{{\mathbb Z}}

\def\nor{\trianglelefteq}
\def\sdp{\rtimes}
\def\nle{\not\le}
\def\iso{\cong}
\def\ord#1{\vert #1 \vert}

\def\gp#1{\langle \, #1 \, \rangle}
\def\phi{\varphi}
\def\Phi{\varPhi}

\def\ZZ{\Bbb Z}

\def\FF{\Bbb F}

\def\O{\mathcal O}

\def\divides{\bigm |}
\def\notdivides{\not\kern 2.2pt\bigm |}

\def\omegabar#1{\frak A_1(#1)}

\title[Cellular structure of classifying spaces]
{The cellular structure of the classifying spaces of finite groups
}

\author{Ram\'on J. Flores and Richard M. Foote}
\thanks{The first author was supported by MEC grant MTM2004-06686.}

\date{\today}
\begin{document}

\begin{abstract}
In this paper we obtain a description of the $B\Z /p$-cellularization
(in the sense of Dror-Farjoun) of the classifying spaces of all finite groups,
for all primes~$p$.

\end{abstract}

\maketitle

\section{Introduction}

Let $A$ be a pointed space. A space $X$ is called $A$-cellular if
it can be constructed as an (iterated) pointed homotopy colimit of
copies of $A$. The concept of $A$-cellularity was developed by
Dror-Farjoun (\cite{Farjoun95}) and Chach\'olski
(\cite{Chacholski96}) with a two-fold goal: to classify the spaces
in cellularity classes, whose properties should depend on those of
the generator $A$, and to develop an $A$-homotopy theory, where
the suspensions of $A$ would play the role that the spheres play
in classical homotopy. In this context these authors constructed
an $A$-cellularization endofunctor $\textbf{CW}_A$ of the category
of pointed spaces that mimics the well-known cellular
approximation.  The functor $\textbf{CW}_A$ turns out to be
augmented and idempotent and is characterized by the facts that
$\map_{*}(A,\textbf{CW}_AX)$ is weakly homotopy equivalent to
$\map_{*}(A,X)$ and every map $A\ra X$ factors through
$\textbf{CW}_AX$ in a unique way, up to homotopy.

These ideas have made a great impact both in Homotopy Theory
(\cite{CDI06}, \cite{CCS2}, \cite{RS08}) and outside of it,
because they can be generalized in a natural way to
any framework where there is a notion of limit (see examples
in \cite{Dwyer06}, \cite{BKI08} or \cite{Kiessling08}). A line of
research of great relevance in the last years, and very related to
our work, is cellular approximation in the category of groups
(\cite{Rodriguez01}, \cite{Ramon}, \cite{DGS07}, \cite{DGS08}).

The present paper culminates a program aimed at understanding the mod
$p$ homotopy theory of the classifying spaces of finite groups $G$
using tools from $\textrm{B}\Z /p$-homotopy theory, where $p$ is any prime.
The main underlying idea in our study is the relationship between the
$\textrm{B}\Z /p$-cellular structure of $BG$ and certain strongly
closed subgroups of $G$. Recall that given a finite group $G$,
a subgroup $R$ of some Sylow $p$-subgroup $S$ of $G$ is
said to be \emph{strongly closed} in $S$ with respect to $G$
if whenever $x\in R$ and $g\in G$ are such that
$gxg^{-1}\in S$, then $gxg^{-1}\in R$.
(All relevant definitions, with some discussion, are collected in
Subsection~\ref{definitions} at the end of this Introduction;
some definitions also appear in the Introduction itself to
maintain the flow of its overview nature.)

This paper gives a complete description of the $B\Z
/p$-cellularization of classifying spaces of all finite groups.
The philosophy behind our work is the following: whenever $X$ is a
space with a notion of $p$-fusion --- and hence strong closure ---
knowledge of the strongly closed subobjects of $X$ is deeply
related (and in some cases, almost equivalent) to the $A$-cellular
structure of $X$, for a certain $p$-torsion space $A$. This
strategy opens up new perspectives for analyzing, from the point
of view of (Co)localization Theory, the $p$-primary structure of a
wide class of homotopy meaningful spaces, such as classifying
spaces of compact Lie groups, $p$-local finite groups, $p$-compact
groups or, more generally, $p$-local compact groups
(\cite{Broto07}).

In the specific case of $G$ finite it became evident from
\cite{Flores07} that the classification of the possible homotopy
types of $\textbf{CW}_{B\Z /p}BG$ would require some kind of
description of the strongly closed subgroups that $G$ can possess.
For $p=2$ this task was completed by the second author in
\cite{Foote97}, while the case of an odd prime was solved in the
separate paper \cite{Flores08}, which is the group-theoretic
underpinning of this work. (The combined classification is rendered
in a slightly simplified form as Theorem~\ref{strongly-closed}
of Section~\ref{examples}.)
The $p$ odd classification was a crucial
ingredient we lacked in order to finish the characterization of
$\mathbf{CW}_{B\Z /p}BG$ for all finite groups $G$ (the other was
the role of the subgroup $\O_A(G)$, see below), solving a problem
that was posed by Dror-Farjoun \cite[3.C]{Farjoun95} in the case
$G=\Z /p^r$, and partially solved in \cite{Ramon} and
\cite{Flores07} (see Section~2 below for an analysis of the
previous cases). The latter paper showed the relationship between
the upper homotopy of $\textbf{CW}_{B\Z /p}BG$ and a specific
strongly closed subgroup $\omegabar S$ of $G$.  More precisely,
for $S$ a Sylow $p$-subgroup of $G$, $\omegabar S$ is the unique
minimal strongly closed subgroup of $S$ that contains all elements
of order $p$ in $S$.
The importance of the subgroup $\omegabar S$ comes from the fact that it
determines a great part of the structure of
$\textrm{Hom} (\Z/p,G)$, and then of $\textrm{map}_*(B\Z /p, BG)$.

We are now in position to explain our main results.
Given a finite group $G$ and any $p$-subgroup $A$ of $G$ then,
as shown in Subsection~\ref{definitions}, there is a unique
subgroup of $G$, denoted by $\O_A(G)$, that is maximal with
respect to the two properties:

\begin{itemize}

\item[(i)]
$\O_A(G)$ is a normal subgroup of $G$, and
\item[(ii)]
$A$ contains a Sylow $p$-subgroup of $\O_A(G)$, i.e., $A \cap \O_A(G) \in Syl_p(\O_A(G))$.
\end{itemize}


\noindent
We are especially interested in the case when $G$ is generated by
its elements of order $p$ and $A = \omegabar S$ for some Sylow
$p$-subgroup $S$ of $G$.  In this situation $A = S$ if and only if
$A \le \O_A(G)$, i.e., if and only if $A$ maps to the trivial
group under the natural projection of $G$ onto $G/\O_A(G)$.
(Section~\ref{examples} and references \cite{Flores07} and
\cite{Flores08} give illuminating examples of the subgroups just
defined.)

Let $BG^{\wedge}_p$ denote Bousfield-Kan $p$-completion of $BG$.
For the definition and properties of this functor, we refer the
reader to \cite{Bousfield72} for a thorough account, and to
\cite[Section 2]{Ramon} for a brief survey. Our first main
result---which is Theorem~\ref{thetheorem} in the paper---is the
following:

\bigskip
\noindent {\bf Theorem A}. \emph{Let $G$ be a finite group
generated by its elements of order $p$, let $A = \omegabar S$ be
the minimal strongly closed subgroup of $S$ containing
all elements of order $p$ in $S$, and assume $A \ne S$.
Let overbars denote passage
from $G$ to the quotient group $G / \O_A(G)$. Then there exists a
fibration sequence
$$
\mathbf{CW}_{B\Z /p}(BG^{\wedge}_p)
\longrightarrow
BG^{\wedge}_p
\longrightarrow
B(N_{\overline{G}}(\overline{A})/ \, \overline{A})^{\wedge}_p.
$$
}

Although it is relatively elementary to prove that for any
strongly closed subgroup $A$ the subgroup $N_G(A)$ controls fusion
of subgroups containing $A$, it does not generally control fusion
inside $A$, or for subsets that intersect $A$ but do not contain
it. The much more delicate problem is to realize that factoring
out the ``correct'' subgroup $\O_A(G)$ combined with an explicit
knowledge of the structure of the quotient $G/\O_A(G)$ is crucial,
in the sense that there is a part of the fusion inside $A$ that
can be ``swept under the rug'' for the topological considerations.
The main result of the classification in \cite{Flores08} is a
determination of the isomorphism type of $\overline G$, and hence
of $N_{\overline G}(\overline A)$, under the hypotheses of
Theorem~A (and more generally as well); roughly speaking, when
$\overline A \ne \overline 1$ the quotient group $\overline G$ is
a direct product of simple groups from certain explicitly listed
families, and in each family the subgroups $\overline A$ and
$N_{\overline G}(\overline A)$ are explicitly determined as well
(a more precise rendition of the classification is given as
Theorem~\ref{strongly-closed}). The full force
of this classification is then invoked to yield the fusion
analysis needed to complete the proofs of Theorem~A and Theorem~B:
namely, that $N_{\overline G}(\overline A)$ always controls
strong fusion in $\overline A$.

In the special case when $G/\O_A(G)$ is just a single simple
group isomorphic to one of the ``obstruction groups''
listed in the conclusions of Theorem~\ref{strongly-closed},
in all but one family the subgroup
$N_{\overline{G}}(\overline{S})$ controls strong fusion in $\overline S$
(and this was the situation for all ``obstructions'' when $p=2$).
The Sylow-normalizer is always
a subgroup of $N_{\overline G}(\overline A)$, however it can be strictly smaller.
In that one exceptional family alone, $N_{\overline G}(\overline S)$
is indeed smaller than
$N_{\overline G}(\overline A)$, and the smaller normalizer does {\it not}
control strong fusion in $\overline A$.
This illustrates that the ultimate classification
for $p$ odd involves more subtle, unavoidable configurations than those that
arose for $p=2$.  These are explicated in more detail in Section~\ref{examples}
and in \cite{Flores08}.

To describe the second main result let $\Omega_1(G)$ denote the
subgroup of $G$ generated by the elements of order $p$ in $G$.
By Propositions~4.1 and 4.3 of \cite{Ramon},
the inclusion $\Omega_1(G)\hookrightarrow G$ induces a
homotopy equivalence
$\mathbf{CW}_{B\Z/p}B\Omega_1(G)\simeq\mathbf{CW}_{B\Z /p}BG$,
so these propositions imply that it is enough to consider the
case in which $G$ is generated by elements of order $p$.
Thus the following result---which is Theorem \ref{recopcel}
herein---combines the information obtained in \cite{Ramon} and
\cite{Flores07} with the present article to give all possible
homotopy structures for $\mathbf{CW}_{B\Z /p}BG$.

\bigskip
\noindent
{\bf Theorem B}. {\it Let $G$ be a finite group
generated by its elements of order $p$, let $S$ be a Sylow $p$-subgroup of $G$,
and let $A = \omegabar S$ be the minimal strongly closed subgroup of
$S$ containing all elements of order $p$ in $S$.
Let overbars denote passage from $G$ to the
quotient group $G / \O_A(G)$. Then the $B\Z/p$-cellularization of
$BG$ has one the following shapes:

\begin{itemize}

\item[(1)]
If $G = S$ is a $p$-group then $BG$ is $B\Z /p$-cellular.

\item[(2)]
If $G$ is not a $p$-group and $A = S$ then $\mathbf{CW}_{B\Z/p}BG$
is the homotopy fiber of the natural map
$BG \ra \prod_{q\neq p}BG^{\wedge}_q$.

\item[(3)]
If $G$ is not a $p$-group and $A \ne S$ then
$\mathbf{CW}_{B\Z/p}BG$ is the homotopy fiber of the map
$BG \ra B(N_{\overline G}(\overline A)/\, \overline A)^{\wedge}_p \times
\prod_{q\neq p}BG^{\wedge}_q$.

\end{itemize}
}

\bigskip

The classification of groups containing a strongly closed
$p$-subgroup in \cite{Flores08} gives a very precise description
of the fiber of the augmentation $\mathbf{CW}_{B\Z /p}BG\ra BG$ in
terms of normalizers of strongly closed $p$-subgroups in the
simple components of $G/\O_A(G)$. In this sense the results here
further improve those of \cite{Flores07}, where this degree of
sharpness was only obtained in the description of some concrete
examples.

The overall organization of the paper is as follows. Section~2
begins by recapitulating previous results around the $B\Z
/p$-cellularization of $BG$. Section~3 compiles some facts about
the structure of mapping spaces from classifying spaces, which are
needed afterwards. Section~4 contains the main results and their
proofs; and Section 5 illustrates the efficacy of our methods by
describing precisely some families of examples. The latter are
very illuminating in the sense that they give an alluring glimpse
of what ``should be'' the $B\Z /p$-cellularization of more general
objects.

\medskip

\subsection{Definitions}
\label{definitions}

Relevant definitions are collected in this subsection,
and in the next subsection we list notation used in this paper.
Readers may also repair to Section~\ref{examples}, independent of the
intervening sections, to plumb some examples that illustrate the concepts.
Throughout the paper $G$ is a finite group, $p$ is a prime and
$S$ is a Sylow $p$-subgroup of $G$.

A subgroup $R$ of $S$ is called
{\it strongly closed} in $G$ if whenever $x \in R$ and $g \in G$ are such that
$gxg^{-1} \in S$, then $gxg^{-1} \in R$ (or equivalently,
whenever $\gp{gxg^{-1},R}$ is a $p$-group, then $gxg^{-1} \in R$,
so this property is independent of $S$).
The {\it $p$-socle} of any group $H$ is the subgroup of $H$ generated by
its elements of order $p$ --- denoted by $\Omega_1(H)$ (here $p$ is fixed
even if $H$ is not a $p$-group).
The subgroup $\omegabar S$ denotes {\it the unique smallest subgroup of $S$
that contains $\Omega_1(S)$ and is strongly closed in $G$}.
Examples of strongly closed subgroups are given in
Section~\ref{examples}.  In particular, if $R = S \cap N$ is a
Sylow $p$-subgroup of a normal subgroup $N$ of $G$, then $R$ is
strongly closed.  This observation indicates that strongly closed
subgroups often, but not always, signal the presence of normal
subgroups.  It also leads naturally to the following ``functor.''

Observe that for $R$ any subgroup of $S$,
if $N_1$ and $N_2$ are normal subgroups
of $G$ with $R \cap N_i \in Syl_p(N_i)$ for both $i=1,2$,
then $R \cap N_1 N_2$ is a Sylow $p$-subgroup of $N_1 N_2$.
Taking the product of all such $N_i$, let $\O_R(G)$ denote the
{\it unique largest normal subgroup $N$ of $G$ for which
$R \cap N \in Syl_p(N)$}.
Thus $\O_R(G)$ is characterized by the two properties given earlier in the
Introduction.  Moreover,
$$
R \text{ is a Sylow $p$-subgroup of $\gp{R^G}$ if and only if } R
\le \O_R (G).
$$
Note that $R \O_R(G) / \O_R(G)$ does not contain the
Sylow $p$-subgroup of any nontrivial normal subgroup of $G/\O_R(G)$;
in other words, $\O_{\overline R}(\overline G) = 1$,
where overbars denote passage to $G/\O_R(G)$.
In the special case when $R$ is strongly closed in $G$,
observe that strong closure passes to all quotient groups,
so when analyzing groups where
$R \nle \O_R(G)$ we often factor out $\O_R(G)$ (and still have a
strongly closed subgroup $\overline R$ of $\overline G$).

For any $R \le G$, we say a subgroup $H$ of $G$ containing $R$
{\it controls fusion in $R$} if any pair of elements of $R$ that are
conjugate in $G$ are also conjugate in $H$.
We likewise say $H$ {\it controls strong fusion} in $R$ if for
every subgroup $P$ of $R$ and $g\in G$ such that $gPg^{-1} \le R$,
there exist $h \in H$ and $c\in C_G(P)$, with $g=hc$.

\subsection{Specific Notation}

Unless explicitly stated otherwise, throughout the paper
we adopt the additional notation that the group $A = \omegabar S$.
(For consistency with cited literature, occasionally the letter $A$
is used as a topological space; but this should not cause
confusion because of the context of such deviations.)
Overbars denote passage to the quotient $G \rightarrow G/\O_A(G)$.
The normalizer of $\overline{A}$ in $\overline{G}$ is denoted by
$\overline{N}$, and the group $\overline{N}/\overline{A}$
will be called~$\Gamma$.

For any space $X$ let $X^{\wedge}_p$ denote the Bousfield-Kan
$p$-completion of $X$. The cofiber of the map $\bigvee B\Z /p\ra
BG$ (where the wedge is defined over all the homotopy classes of
maps from $B\Z /p$ to $BG$) is called $C$. Analogously, the
cofiber of the corresponding map $\bigvee B\Z /p\ra BG^{\wedge}_p$
(the wedge defined over classes of maps $B\Z/p\ra BG^{\wedge}_p$)
is denoted by $D$. Let $P=\mathbf{P}_{\Sigma B\Z /p}C$, the
$\Sigma B\Z /p$-nullification of $C$ (see the remarks following
Theorem~\ref{theorem2-1} for the definition); and it was
checked in the proof of \cite[Proposition 3.2]{Flores07} that
$\mathbf{P}_{\Sigma B\Z /p}D$ has the homotopy type of
$P^{\wedge}_p$, so we will refer to this object with this name.

\medskip
\noindent {\bf Acknowledgements}. We warmly thank Carles Broto,
Bob Oliver and J\'er\^ome Scherer for helpful discussions. We
thank the referee for many helpful and insightful suggestions that
greatly improved the presentation.


\section{Previous results}

Before undertaking the complete description of $\mathbf{CW}_{B\Z /p }BG$
we describe what is known so far about this problem.
As we said in the Introduction, the starting point was the
computation done by Dror-Farjoun in \cite[3.C]{Farjoun95}, where
he establishes that the $B\Z /p$-cellularization of the
classifying space of a finite cyclic $p$-group has the homotopy
type of $B\Z /p$.

Subsequently Rodr{\'\i}guez-Scherer investigated in
\cite{Rodriguez01} the $M(\Z /p,1)$-cellularization, where $M(\Z
/p,1)$ denotes the corresponding Moore space for $\Z /p$. When the
target is $BG$, this can be considered a precursor to our
study because $M(\Z /p, 1)$ can be described as the 2-skeleton of
$B\Z /p$. In their description the authors use the concept of
cellularization in the category of groups (developed afterwards in
\cite{Farjoun07}). Their work in this subject allows one to prove, in
particular, that the $B\Z /p$-cellularization of the classifying
space of a $p$-group is the same as that of its $p$-socle; as the
latter is $B\Z /p$-cellular in this case (\cite[Proposition
4.14]{Ramon}), one obtains that $\mathbf{CW}_{B\Z /p}BG\simeq
B\Omega_1(G)$ if $G = S$ is a finite $p$-group.

The aforementioned is proved using a characterization of
the cellularization discovered by Chach\'olski, that is perhaps the
most useful tool available to attack these kind of problems.
Because of its importance and ubiquity in our context
we reproduce it here:


\begin{theorem} \cite[20.3]{Chacholski96} Let $A$ and $X$ be pointed
spaces, and let $C$ be the homotopy cofiber of the map
$\bigvee_{[A,X]_*} A\ra X$, defined as evaluation over all the
homotopy classes of maps $A\ra X$. Then $\mathbf{CW}_{A}X$ has the
homotopy type of the homotopy fiber of the composite $X\ra
C\ra\mathbf{P}_{\Sigma A }C.$ \label{Chacholski}
\label{theorem2-1}
\end{theorem}


Here $\mathbf{P}$ denotes the nullification functor, first defined
by A.K. Bousfield in \cite{Bousfield94}. Recall that given spaces
$A$ and $X$, $X$ is called $A$-\emph{null} if the natural
inclusion $X\hookrightarrow\textrm{map}(A,X)$ is a weak
equivalence. In this way one defines a functor
$\mathbf{P}_A:\mathbf{Spaces}\ra\mathbf{Spaces}$, coaugmented and
idempotent, such that $\mathbf{P}_AX$ is $A$-null for every $X$,
and such that for every $A$-null space $Y$ the coaugmentation
induces a weak equivalence
$\textrm{map}(\mathbf{P}_AX,Y)\simeq\textrm{map}(X,Y)$. This
functor can also be defined in the pointed category, and its
main properties can be found in \cite{Farjoun95} and
\cite{Chacholski96}.

In our case the role of $A$ and $X$ in Theorem \ref{Chacholski}
will be played by $B\Z /p$ and $BG$, respectively. If $C$ is the
corresponding cofiber, from now on we shall denote the $\Sigma B\Z
/p$-nullification of $C$ by $P$. As a consequence of Theorem
\ref{Chacholski}, describing $\mathbf{CW}_{B\Z /p}BG$ is
equivalent to describing $P$, which is in general a more
accessible problem. Next we list some properties of the space $P$
which will be useful in the remaining of our note.

\begin{proposition}

Let $P$ be as described in the previous paragraph, with $G$ a
finite group generated elements by order $p$. Then the following
hold:

\begin{enumerate}

\item The spaces $P$ and $P^{\wedge}_p$ are simply-connected.

\item $P^{\wedge}_p$ is $H\mathbb{Z}/p$-local.

\item The loops of $P$ and $P^{\wedge}_p$ are $B\Z /p$-null
spaces.

\item $P$ admits a splitting $(\prod_q
BG^\wedge_q) \times P^{\wedge}_p$, where the product is taken over
all primes $q \neq p$; moreover $P^{\wedge}_p$ coincides with the
base of the fibration sequence of Theorem \ref{Chacholski} when
the total space is $BG^{\wedge}_p$.

\end{enumerate}

\begin{proof}
(1) $C$ is by definition a homotopy push-out ${*}\la \bigvee B\Z
/p \ra BG$, and the generation hypothesis on $G$ implies that
the right map is surjective at the level of fundamental groups.
Hence by Seifert-Van Kampen $C$ is simply-connected, and then $P$
is so by \cite[Proposition 6.9]{Dwyer96}. In particular, by
\cite[VII,3.2]{Bousfield72}, $P^{\wedge}_p$ is also 1-connected.

(2) According to \cite[VI.5.1]{Bousfield72} and
\cite[Proposition 4.3]{Bousfield75}, the homotopy groups of
$P^{\wedge}_p$ are $H\Z /p$-local in the sense of Bousfield, and
hence $P^{\wedge}_p$ is $H\Z /p$-local (as a space) by
\cite[Theorem 5.5]{Bousfield75}.

(3) By adjunction, $\textrm{map}_{*}(B\Z /p,\Omega P)$ is weakly
equivalent to $\textrm{map}_{*}(\Sigma B\Z /p,P)$, which is weakly
trivial, and then $\Omega P$ is $B\Z /p$-null. The fibre lemma
\cite[II.5.1]{Bousfield72} applied to the loop fibration of $P$
implies that $\Omega (P^{\wedge}_p) \simeq (\Omega P)^{\wedge}_p$,
and the conclusion follows from \cite[Lemma 9.9]{Miller}.

(4) See Proposition 2.1 and the proof of Proposition 3.2 in \cite{Flores08}.

\end{proof}
\label{cofibre}
\end{proposition}

Note that the additional hypothesis on the generation of $G$
causes no restriction from the point of view of cellularization,
as for every finite $G$ there is a homotopy equivalence
$\mathbf{CW}_{B\Z /p}B\Omega_1(G)\simeq\mathbf{CW}_{B\Z /p}BG$
induced by the natural inclusion (\cite[Proposition 4.1]{Ramon}).
So, one should have in mind that the knowledge of
$\mathbf{CW}_{B\Z /p}BG$ for $G$ generated by order $p$ elements
implies automatically the knowledge of $\mathbf{CW}_{B\Z /p}BG$
for every $G$.

According to \cite[3.2]{Flores07}, the last statement of our
previous Lemma implies that we are actually studying the $B\Z
/p$-cellularization of $BG^{\wedge}_p$. More precisely, if $G$ is
generated by order $p$ elements, the equivalence $\mathbf{CW}_{B\Z
/p}(BG^{\wedge}_p)\simeq ( \mathbf{CW}_{B\Z /p}BG)^{\wedge}_p$ is
proved in \cite[Proposition 3.2]{Flores07}. This is a very
peculiar property for spaces which, in general, cannot be
decomposed via an arithmetic square.

In the philosophy of \cite{Broto032}, the homotopy theory of
$BG^{\wedge}_p$ is codified in the $p$-fusion data of $G$. From
this point of view it can be observed that the structure of
$P^{\wedge}_p$ strongly depends on the minimal strongly closed
$p$-subgroup $\omegabar S$ of $S$ that contains the $p$-socle of
$S$ (called $Cl\textrm{ }S$ in \cite{Flores07}). In particular, it
is a consequence of the Puppe sequence and the definition of
nullification (see \cite[Proposition 4.2]{Flores07} for details)
that if $\omegabar S=S$ then $P^{\wedge}_p$ is trivial. This
leads one to consider the case in which $\omegabar S$ is strictly
contained in $S$.

In \cite{Flores07} the latter case is studied under the additional
assumption that $N_G(S)$ controls (strong) $G$-fusion in $S$.
Since
$\omegabar S$ is normal in $N_G(S)$, \cite{Flores07} shows that
$P^{\wedge}_p$ is homotopy equivalent to the $p$-completion of
$B(N_G(S)/\omegabar S)$; this also shows, roughly speaking, that
the structure of the mapping space $\textrm{map}_*(B\Z
/p,BG^{\wedge}_p)$ depends heavily on $\omegabar S$. This result
is used, in particular, to compute the $B\Z /2$-cellularization of
classifying spaces of simple groups (relying on Theorem 5.6
there).

The explicit description of $P^{\wedge}_p$ in the remaining cases,
which lead in turn to a complete knowledge of the structure of
the  $B\Z /p$-cellularization of $BG$ for every finite group $G$, is
given in Section~4. Before, however, we need to deal with some
mapping spaces which play a crucial role in the description.

\section{Maps from classifying spaces and Zabrodsky Lemma}

It was already evident in our previous work \cite{Flores07} that a
complete knowledge of the homotopical structure of
$\mathbf{CW}_{B\Z /p}BG$ demands, for some spaces $X$, a certain
control of the behavior of the functor $\textrm{map}_*(-,X)$ when
applied to fiber sequences of classifying spaces of finite groups.
This is done in part by the classical Zabrodsky Lemma, which we recall
here for the reader's convenience, from a recent version of W.~Dwyer. We
also present in this section a specific version of this statement,
which will be useful in the proof of our key result (Theorem
\ref{thetheorem}).

Let $F\ra E\ra B$ be a fibre sequence, with $B$ connected. Denote
by $\textrm{map}(F,X)_{[F]}$ the space of maps $F\ra X$ which are
null-homotopic, by $\textrm{map}(E,X)_{[F]}$ the space of maps
$E\ra X$ which become null-homotopic when restricted to $F$, and
by $[E,X]_{[F]}$ the group of components of
$\textrm{map}(E,X)_{[F]}$, which is identified with the homotopy
classes of maps $E\ra X$ which restrict trivially to $F$; observe
in particular that all these concepts have their pointed
counterparts. Then we have the following:

\begin{theorem}[\cite{Dwyer96}, Proposition 3.5]

If $X$ is a space and the inclusion $X\ra\emph{map}(F,X)_{[F]}$ is an
equivalence, then the restriction
$\emph{map}(B,X)\ra\emph{map}(E,X)_{[F]}$ is an equivalence too.

\label{DwyerZabrodsky}
\end{theorem}

Before introducing our version of this result, which is focused on
the case in which the fibre sequence is derived from an extension
of groups, we recall some well known properties of the
mapping spaces and the $\Sigma B\Z /p$-null spaces. This is
undertaken in the next two lemmas.

\begin{lemma}

Let $A$ and $X$ be spaces with $X$ simply-connected. Then there
is a bijection $[A,X]_*\ra [A,X]$ between the set of components of
the respective mapping spaces $\emph{map}_{*}(A,X)$ and
$\emph{map}(A,X)$.

\label{pointedunpointed}
\end{lemma}

\begin{proof}

It is enough to consider the long exact homotopy sequence of the
classical fibre sequence $$\textrm{map}_{*}(A,X)\ra
\textrm{map}(A,X)\ra X,$$ and to take account of the fact that
$\pi_0\textrm{map}_{*}(A,X)=[A,X]_*$ and
$\pi_0\textrm{map}(A,X)=[A,X].$
\end{proof}

\begin{lemma}

Let $p$ be a prime. If $G$ is a compact Lie group and $X$ is a
simply-connected $p$-complete, $\Sigma B\Z /p$-null space, then
$X$ is $\Sigma BG$-null.

\label{sigmaBGnull}
\end{lemma}

\begin{proof}

Recall that by definition of nullification, a space $X$ is $\Sigma
A$-null for a certain space $A$ if and only if $\Omega X$ is
$A$-null. Then $\Omega X$ is $p$-complete and $B\Z /p$-null, so
according to \cite[Theorem 1.2]{Dwyer96} it is $BG$-null. Hence
$X$ is $\Sigma BG$-null.
\end{proof}

We may now state our specific version of Zabrodsky's Lemma:

\begin{proposition}

Let $G_1\ra G_2\stackrel{h}{\ra }G_3$ be an extension of finite
groups, $BG_1\ra BG_2\ra BG_3$ the corresponding fibre sequence of
classifying spaces. If $X$ is a simply-connected, $p$-complete
$\Sigma B\Z/p$-null space for a prime $p$, $f$ induces a bijection
$[BG_3,X]\simeq [BG_2,X]_{[BG_1]}$, which remains true in the
pointed category. \label{ZabrodskySpecific}
\label{zabrodskypointed}
\end{proposition}

\begin{proof}

We follow the proof of \cite[Lemma 2.3]{CCS2}. By the previous
Lemma \ref{sigmaBGnull}, $X$ is $\Sigma BG_1$-null, so by
definition $\textrm{map}_*(\Sigma BG_1,X)$ is contractible, and in
particular the component of the constant in
$\textrm{map}_*(BG_1,X)$ is so. Then, the restriction of the fibre
sequence of the proof of \ref{pointedunpointed} (with $A=BG_1$) to
the components of constant maps gives that the evaluation
$\textrm{map}(BG_1,X)_{[BG_1]}\ra X$ is an equivalence. Now
Theorem \ref{DwyerZabrodsky} implies that the map
$\textrm{map}(BG_3,X)\ra \textrm{map}(BG_2,X)_{[BG_1]}$ induced by
$BG_2\stackrel{Bh}{\ra }BG_3$ is also an equivalence, and
therefore $[BG_3,X]\simeq [BG_2,X]_{[BG_1]}$. As $X$ is
simply-connected, Lemma \ref{pointedunpointed} ensures that the
bijection remains true for classes of pointed maps, so we are
done.
\end{proof}

This proposition can be read, in particular, as an extension
property:

\begin{corollary}

In the conditions of the previous proposition, if $f:BG_2\ra X$ is
a map in $\emph{map}(BG_2,X)_{[BG_1]}$, then $f$ factors through a
map $g:BG_3\ra X$ that is unique up to unpointed homotopy.
Moreover, if we assume $f$ to be pointed, then the factorization
is unique up to pointed homotopy.

\label{corozabrodsky}
\end{corollary}

\begin{remark}

Our usual candidates for $X$ will be $p$-completions of
classifying spaces of $p$-perfect finite groups, which are
simply-connected \cite[Proposition II.5]{Bousfield72},
$p$-complete \cite[Proposition II.5]{Bousfield72} and $\Sigma B\Z
/p$-null \cite[Lemma 9.9 ]{Miller}, the latter because $BG$ is
always so. We will make regular use of this fact.

\label{Mill}
\end{remark}

We include a group-theoretic lemma that will be needed
in the proof of Proposition~\ref{CCS2-2-4} following,
and also in the description of the structure of $A = \omegabar S$
in Section~4.

\begin{lemma}
\label{subgroupindexp} Let $T$ be a nontrivial $p$-group, let $\{
x_1, \dots, x_k \}$ be any minimal set of generators of $T$, and
let $H$ be the normal subgroup generated by all $T$-conjugates of
$x_1,\dots,x_{k-1}$. Then $H$ is a proper subgroup of $T$ with $T
= H\gp{x_k}$, and $T/H$ is isomorphic to a (cyclic) quotient group
of $\gp{x_k}$. In particular, if all $x_i$ have order $p$, then $H
= \Omega_1(H)$ and $T/H \iso \ZZ/p$.
\end{lemma}

\begin{proof}

Let $x_1,\dots,x_k$ be a minimal set of generators of $T$, let
$H_0 = \gp{x_1,\dots,x_{k-1}}$, and let $H$ be the normal closure
of $H_0$ in $T$. By minimality of $k$ we have $H_0 < T$. Then
$H_0$ lies in some maximal subgroup, $M$, of $T$, and by basic
$p$-group theory $M \nor T$; thus we have $H \le M < T$ too. Since
$T = \gp{H_0,x_k}$ we get $T = H\gp{x_k}$; and so $T/H$ is
isomorphic to a quotient of $\gp{x_k}$.  This proves the first
assertion.

If all $x_i$ have order $p$, then since $H$ is generated by
conjugates of these, $H = \Omega_1(H)$, and the second
assertion of the lemma now follows from the first.
\end{proof}

We finish the section with an adaptation of Proposition 2.4 in
\cite{CCS2}, which will be necessary when dealing with maps whose
source is a strongly closed subgroup.

\begin{proposition}
\label{CCS2-2-4}

Let $P\stackrel{i}{\ra} Q\stackrel{p}{\ra} G$ be an extension of
finite groups, $X$ a simply-connected $\Sigma B\Z /p$-null
$p$-complete space, and $f:BQ\ra X$ a map that restricts
trivially to $BP$. Assume that $G$ is a $p$-group generated by a
collection $C$ of elements such that, for every $x\in C$, there
exists $y\in Q$ such that $p(y)=x$ and the restriction of $f$ to
$B\langle y\rangle$ is null-homotopic. Then $f$ is null-homotopic.

\label{extensionCCS}
\end{proposition}

\begin{proof}

According to the previous corollary, there exists an extension
$\hat{f}:BG\ra X$ of $f$ such that $f$ is null-homotopic if and
only if $\hat{f}$ is. We will check the latter by induction on the
order of the group $G$. If $G=\langle x\rangle$ with $x$ an
element of order $p$, then the composite $B\langle y\rangle\ra
B\langle x\rangle\ra X$ is null-homotopic, and $B\langle
x\rangle\ra X$ is so by applying Corollary \ref{corozabrodsky} to
the fibration $BK\ra B\langle y\rangle\ra B\langle x\rangle$, with
$K$ the kernel of $\langle y\rangle\ra \langle x\rangle$. So we
are done in this case.

For the general case, apply Lemma~\ref{subgroupindexp}
to $T = G$ and the set $C =\{x_1,\ldots ,x_k\}$,
which we may assume is a minimal set of generators.
We obtain a proper normal subgroup $H$ of $G$
and an extension $H\ra G\ra \Z /{p^r}$ with
$r\geq 1$, where the quotient is generated by the image of $x = x_k$.
Let $Q'$ be the pull-back of the extension of the statement along the
inclusion $H\ra G$, let $Q'\ra Q$ be the corresponding induced
homomorphism, and let $h:BQ'\ra BQ$ be the induced map between
classifying spaces. Now $H$ and the extension $P\ra Q'\ra H$
satisfy the assumptions in the statement, and moreover $|H|<|G|$
so by induction the composition $f\circ h$ is trivial,
and the restriction of $\hat{f}$ to $BH$ is so. Now consider the
diagram:

$$
\xymatrix{
B(\langle x\rangle\cap X) \ar[r] \ar[d] & BH \ar[d]
\ar[dr]^{{x}} & \\
B(\langle x\rangle) \ar[r] \ar[d] & BG \ar[d] \ar[r]^{\hat{f}} & X
\\
B\Z /{p^r} \ar[r]^{Id} & B\Z /{p^r} \ar[ur]_{f'} & }
$$

By Corollary \ref{corozabrodsky} applied to the second vertical
fibre sequence, we only need to see that $f'$ is trivial, and
again by the same result applied now to the first, we obtain that
$f'$ is trivial because $\hat{f}$ is trivial too. So $\hat{f}$ is
null-homotopic and we are done.
\end{proof}


\section{Describing the $B\Z /p$-cellularization of $BG$}

Unless otherwise stated, throughout this section $G$ will always
be a finite group generated by its elements of order $p$. As said above,
this is a technical reduction which does not affect the generality
of our main result.  We continue to use notation specified
at the end of Section~1.

In \cite{Flores07} it was already anticipated that a complete
description of the $B\Z /p$-cellularization of classifying spaces
of finite groups would depend on a structure theorem for groups
that contain a non-trivial strongly closed $p$-subgroup that is
not a Sylow $p$-subgroup; at that time such a classification was
only known for $p=2$ (\cite{Foote97}). But even in this case there
were examples of groups such that $\omegabar S \ne S$ and $N_G(S)$
does not control fusion in $S$
--- in other words, groups that were beyond the scope of
\cite{Flores07}.

The key step missing in that paper was the role of the subgroup
$\O_A(G)$, whose importance was already evident, in an
independent group-theoretic context, in \cite{Foote97}. First,
$\O_A(G)$ is by definition a subgroup of $G$, so in the case
$A=\omegabar{S}$ it is likely that it controls a ``part'' of the
structure of $\mathbf{CW}_{B\Z /p}BG$. Moreover, $\O_A(G)$ is
normal in $G$, so one might expect a strong relationship between
the cofiber of Chach\'olski fibration in Theorem \ref{Chacholski}
for $BG^{\wedge}_p$ and for $B(G/\O_A(G))^{\wedge}_p$. The
significance of $\O_A(G)$ is also suggested by the fact that
passing to the quotient $G/\O_A(G)$ gives a considerable
simplification in the $p$-fusion structure of this quotient group
in the sense that
--- as we shall see in Section~\ref{examples} --- there are
normalizers (of subgroups/elements) that do not control fusion in
$G$ with images that do so in the quotient. These ideas, combined
with the explication of the role of $\O_A(G)$ in the
classification result \cite[Theorem 1.2]{Flores08}, allow us to
prove the main theorem, Theorem~\ref{thetheorem}, which covers all
extant cases for $\mathbf{CW}_{B\Z /p}BG$, subsuming all by a
uniform treatment.

\medskip

At this point it is worth recalling from \cite{Flores07} the
inductive construction of $A$. For $S$ a Sylow $p$-subgroup of $G$
define $Cl_0(S)=\Omega_1(S)$, the subgroup of $S$ generated by its
elements of order $p$. Let $Cl_{i+1}(S)$ be built from
$Cl_{i}(S)$ as the group generated by all elements
$gxg^{-1}\in S$ with $g\in G$
and $x\in Cl_{i}(S)$. Then $A$ is the (finite) union
$\cup_iCl_{i}(S)$, and hence is the unique minimal strongly
closed subgroup of $S$ containing all elements of order $p$ in $S$.
We have the following:

\begin{lemma}

Let $p$ be a prime, $X$ a simply-connected $p$-complete, $\Sigma
B\Z /p$-null space, and let $\gamma:BA\ra X$ be a map that
restricts trivially to every map $B\Z /p\ra BA$. Then $\gamma$ is
also null-homotopic.

\label{TrivialRestriction}
\end{lemma}

\begin{proof}

Suppose first that $A$ is generated by elements of order $p$. We
use induction on the order of $A$. If $|A|=p$, the result is
clear. If not, suppose the result is true for every group $A_1$ such
that $\Omega_1(A_1) = A_1$ and $\ord{A_1} < \ord A$.
As $A$ is a $p$-group generated
by elements of order $p$, Lemma~\ref{subgroupindexp} gives an
extension $H\ra A\ra A/H$ such that $H=\Omega_1(H)$, $H \nor A$ and
$A/H=\Z /p$. Then, applying the induction hypothesis and
Corollary \ref{corozabrodsky} to the fibre sequence $BH \ra BA\ra
B(A/H)$, we obtain that, up to homotopy, there is an extension
$$
\xymatrix{
BA \ar[d] \ar[r] & X \\
B(A/H) \ar@{-->}[ur] & \\
}
$$
and the horizontal map is trivial if and only if the diagonal
one is so. But the latter holds by the induction hypothesis applied
now to $A/H=\Z /p$, so we are done with this case.

Consider now the general case, with $A$ not necessarily generated
by elements of order $p$. We must see that $\gamma$ is trivial
when restricting to $Cl_{i}(S)$ for every $i$. The statement is
clear for $i=0$, because this subgroup is generated by elements of
order $p$. Assume by induction that the restriction to
$BCl_{i}(S)$ for a certain fixed $i$ is null-homotopic, and
consider the fibre sequence $BCl_{i}(S)\ra BCl_{i+1}(S)\ra
B(Cl_{i+1}(S)/Cl_{i}(S))$, which is defined because
$Cl_{i}(S)\trianglelefteq Cl_{i+1}(S)$. By induction
the composite  $BCl_{i}(S)\ra BCl_{i+1}(S)\ra
A\stackrel{\gamma}{\ra} X$ is inessential. Now Corollary
\ref{corozabrodsky} applied to the previous fibre sequence gives
that the map $BCl_{i+1}(S)\ra X$ factors through a map
$s:B(Cl_{i+1}(S)/ Cl_{i}(S))\ra X$, and it is null-homotopic if
and only if $s$ is. But generators of the quotient group
$Cl_{i+1}(S)/Cl_{i}(S)$ are classes of conjugates of elements of
$Cl_{i}(S)$ on which $\gamma$ restricts trivially, so $s$ is
homotopic to the constant when restricted to these generators.
Then the previous fibre sequence satisfies the conditions of
Lemma~\ref{extensionCCS}, and therefore $s$ is null-homotopic,
as needed to complete the induction.
\end{proof}


In the proof of the Theorem \ref{thetheorem}, we will also need a
group-theoretic property of the quotient group
$\Gamma=\overline{N}/\overline{A}$.

\begin{lemma}
\label{pperfect}

The group $\Gamma$ is $p$-perfect, i.e., has no normal subgroup of
index $p$.

\end{lemma}

\begin{proof}

By contradiction, assume that $\Gamma$ is not $p$-perfect, and let
$\psi:\overline{N}\ra \Z /p$ be a surjection that sends
$\overline{A}$ to the identity. As $\overline{N}$ controls fusion
in $\overline{G}$ and the inclusion $\overline{N}<\overline{G}$ is
a mod $p$ homology isomorphism, $\psi$ factors through a
homomorphism on $\overline{G}$. Since $G$ is generated by its
elements of order $p$ and $A$ contains all elements of order $p$
in some Sylow $p$-subgroup of $G$, no proper normal subgroup of
$G$ contains $A$, and then no proper normal subgroup of
$\overline{G}$ contains $\overline{A}$. This is a contradiction,
and we are done.
\end{proof}

%

Alternatively, Theorem B of \cite{Gol75} says that for any
strongly closed $A$ we have $(G' \cap S)A = (N_G(A)' \cap S)A$.
Thus if $G$ has no normal subgroup of index $p$ containing $A$,
neither does $N_G(A)/A$. Corollary~1.5 in \cite{Flores08}
therefore also shows $\Gamma$ is $p$-perfect since each $L_i$ is
simple (see also Corollary~\ref{control-fusion}).

\medskip

Now we are in a position to prove the principal result of this
section. Recall that $G$ is generated by elements of order $p$,
and $P$ denotes the $\Sigma B\Z /p$-nullification of the homotopy
cofibre of the map $\bigvee B\Z /p\ra BG$, where the wedge is
taken over all the homotopy classes of maps $B\Z /p\ra
BG^{\wedge}_p$. We call $D$ this homotopy cofibre. Again overbars
denote passage to the quotient $G \rightarrow G/\O_A(G)$, $A$ is
the unique smallest subgroup of the $p$-Sylow $S<G$ that contains
$\Omega_1(S)$ and is strongly closed in $G$, and $\overline{N}$ is
the normalizer of $\overline{A}$ in $\overline{G}$.

\begin{theorem}

In the previous notation, $P^{\wedge}_p$ is homotopy equivalent
to the $p$-completion of the classifying space of
$\overline{N}/\overline{A}$.

\label{thetheorem}
\end{theorem}

\begin{proof}

We denote by $h:BG^{\wedge}_p\ra D$ the natural map, and by
$\eta:D\ra P^{\wedge}_p$ the canonical coaugmentation. Moreover,
if $A_1<A_2$, we will call $i_{A_1,A_2}$ the group inclusion
\text{$A_1\hookrightarrow A_2$}. Unless it is stated otherwise, we
assume throughout this proof that we work in the pointed category.
Observe that when the target space $X$ is $p$-complete, $\Sigma
B\Z /p$-null and simply-connected, we can use the ``pointed
version" of the Zabrodsky Lemma given by Proposition
\ref{zabrodskypointed} and Corollary \ref{corozabrodsky}.

We claim there are maps $B\Gamma^{\wedge}_p\ra P^{\wedge}_p$ and
$P^{\wedge}_p\ra B\Gamma^{\wedge}_p$ that are homotopy inverses
to one another.
First we define $P^{\wedge}_p\xrightarrow{g} B\Gamma^{\wedge}_p$.
Recall that, as $\overline{N}$ controls $\overline G$-fusion in
$\overline S$, the inclusion $\overline{N}\hookrightarrow
\overline G$ induces a homotopy equivalence
$B\overline{N}^{\wedge}_p\stackrel{(Bi_{\overline{N},\overline{G}})^{\wedge}_p}{\simeq}
B\overline{G}^{\wedge}_p$ (see for example \cite[Proposition
2.1]{MP98}) . Now consider the diagram
\begin{equation} \label{diag1} \xymatrix{ \bigvee B\Z /p
\ar[r]^{v} & BG^{\wedge}_p \ar[d]_{B\pi^{\wedge}_p} \ar[r]^{h} & D
\ar@{-->}[lddd]^{g'}
\ar[r]^{\eta} & P^{\wedge}_p \ar@{-->}[llddd]^{g} \\
& B\overline{G}^{\wedge}_p & & \\
& B\overline{N}^{\wedge}_p
\ar[u]^{(Bi_{\overline{N},\overline{G}})^{\wedge}_p}_{\simeq} \ar[d]_{Bp^{\wedge}_p} & & \\
& B\Gamma^{\wedge}_p & & } \end{equation} and call $\alpha$ the
composite of all vertical maps.

According to the definition of $\overline{G}$ and $A$, the
composite $B\Z /p\ra BG^{\wedge}_p\xrightarrow{\alpha}
B\Gamma^{\wedge}_p$ is inessential for every map $B\Z /p\ra
BG^{\wedge}_p$. This implies that the composite $\alpha\circ v$ is
so, and hence there exists the lifting $g'$. Now Remark \ref{Mill}
and the universal properties of $p$-completion and $\Sigma B\Z
/p$-nullification imply that $g'$ also lifts to $g$, and that is
the map we were looking for.

To construct $f:B\Gamma^{\wedge}_p\ra P^{\wedge}_p$ consider now
the composite $BA\xrightarrow{(Bi_{A,G})^{\wedge}_p} BG^{\wedge}_p
\xrightarrow{\eta\circ h } P^{\wedge}_p$. As the induced
homomorphism of fundamental groups is trivial by construction when
restricted to every element of order $p$ of $A$, Lemma
\ref{TrivialRestriction} implies that the composite is
null-homotopic. In particular, it is also inessential when
precomposing with the map $B(\O_A(G)\cap
A)\xrightarrow{Bi_{\O_A(G)\cap A,A}} BA$. As $P^{\wedge}_p$ is
$p$-complete (by \cite[3.2]{Flores07}), $\O_A(G)\cap A$ is
$p$-Sylow in $\O_A(G)$ by definition of $\O_A(G)$, and the
hypothesis of \cite[Theorem 1.4]{Dwyer96} hold by Proposition~\ref{cofibre},
we obtain that the composite
$B\O_A(G)\xrightarrow{(Bi_{\O_A(G),G})^{\wedge}_p}
BG^{\wedge}_p\longrightarrow P^{\wedge}_p$ is again homotopically
trivial. Then by Proposition \ref{ZabrodskySpecific} applied to
the vertical fibre sequence below, there exists a lifting $f'$
$$
\xymatrix{ B\O_A(G) \ar[d]_{Bi_{\O_A(G),G}} & \\
BG \ar[d]_{B\pi} \ar[r]^{\eta\circ h\circ (-)^{\wedge}_p} & P^{\wedge}_p \\
B\overline{G} \ar@{-->}[ur]_{f'} & }
$$ making the triangle homotopy commutative; here $(-)^{\wedge}_p$ denotes the $p$-completion $BG\ra BG^{\wedge}_p$.
Now we have another map $B\overline{N}\ra P^{\wedge}_p$ given by
the composite
$$
B\overline{N}\longrightarrow B\overline{G}\xrightarrow{f'}
P^{\wedge}_p,
$$
where we have used the fact that $P^{\wedge}_p$ is $p$-complete.

The next diagram is clearly commutative by construction:
$$
\xymatrix{ BA \ar[d]_{B\pi_{|A}} \ar[rr]^{Bi_{A,G}} & & BG
\ar[d]_{B\pi} \ar[dr] \\
B\overline{A} \ar[r]^{Bi_{\overline{A},\overline{N}}} &
B\overline{N} \ar[r]^{Bi_{\overline{N},\overline{G}}}&
B\overline{G} \ar[r] & P^{\wedge}_p.}
$$

Note that the composite
$BA\xrightarrow{Bi_{A,G}}BG\xrightarrow{\eta\circ h\circ
(-)^{\wedge}_p} P^{\wedge}_p$ is homotopically trivial, and hence
the composite $BA\longrightarrow
B\overline{A}\xrightarrow{Bi_{\overline{A},\overline{N}}}
B\overline{G}\stackrel{f'}{\longrightarrow}P^{\wedge}_p$ is too.
By Proposition \ref{ZabrodskySpecific} the composite
$B\overline{A}\xrightarrow{Bi_{\overline{A},\overline{N}}}
B\overline{N}\xrightarrow{Bi_{\overline{N},\overline{G}}}
B\overline{G}\stackrel{f'}{\longrightarrow}P^{\wedge}_p$ is also
homotopically trivial, and thus again by Corollary
\ref{corozabrodsky} applied to the extension
$\overline{A}\ra\overline{N}\ra\Gamma$, the map factors through
$B\Gamma$. As $P^{\wedge}_p$ is $p$-complete, we obtain a map $f$
that fits in the following commutative diagram:
$$
\xymatrix{ B\overline{N}^{\wedge}_p \ar[d]
\ar[r]^{(f')^{\wedge}_p} & P^{\wedge}_p
\\
B\Gamma^{\wedge}_p. \ar[ur]_{f} & }
$$
This is the map $f$ that we wanted. As $B\overline{N}^{\wedge}_p$
is homotopy equivalent to $B\overline{G}^{\wedge}_p$, we have in
particular another commutative diagram, where $\alpha$ is defined
in the obvious way by composing a homotopy inverse of the previous
equivalence with the projection $B\pi$:

\begin{equation} \label{diag6} \xymatrix{ BG^{\wedge}_p
\ar[d]_{\alpha} \ar[r]^{\eta\circ h} & P^{\wedge}_p
\\
B\Gamma^{\wedge}_p. \ar[ur]_{f} & }
\end{equation}

It remains to prove that $f\circ
g\simeq\textrm{Id}_{P^{\wedge}_p}$ and $g\circ f\simeq
\textrm{Id}_{B\Gamma^{\wedge}_p}.$
In the first case, the universal property of the nullification
functor implies that it is enough to prove that $f\circ g'$ is
homotopic to the coaugmentation $\eta$. Now we only need to prove
that $f\circ g'\simeq\eta$; and moreover, as $P^{\wedge}_p$ is
also $\Sigma B\Z /p$-null, we can use the Puppe sequence of the
cofibre sequence $BG^{\wedge}_p\xrightarrow{h} D \ra \vee\Sigma
B\Z /p$ to establish that $f\circ g'\simeq\eta$ if and only if
$f\circ g'\circ h\simeq\eta\circ h$. According to
diagram~\ref{diag1}, $g'\circ h$ is homotopic to $\alpha$, and by
diagram~\ref{diag6}, $f\circ\alpha\simeq\eta\circ h$, so we are
done.

Let us see that $g\circ f$ is homotopic to the identity of
$B\Gamma^{\wedge}_p$. We have the following commutative diagram:
$$
\xymatrix{ BG \ar[r] \ar[d]_{B\pi} & BG^{\wedge}_p \ar[d]^{B\pi^{\wedge}_p} & \\
B\overline{G} \ar [r] &  B\overline{G}^{\wedge}_p \ar@/_/[d]_{j} &
\\
B\overline{N} \ar[u]^{Bi_{\overline{N},\overline{G}}}
\ar[d]_{B\rho} \ar[r] & B\overline{N}^{\wedge}_p
\ar@/_/[u]_{Bi_{\overline{N},\overline{G}^{\wedge}_p}} \ar[d]^{B\rho^{\wedge}_p} & \\
B\Gamma \ar[r] & B\Gamma^{\wedge}_p \ar@/^1pc/[dr]^{Id} \ar@/_1pc/[dr]_{g\circ f} & \\
 & & B\Gamma^{\wedge}_p. }
$$
Here the horizontal arrows represent the $p$-completion map
$(-)^{\wedge}_p$, the homomorphism $\rho:\overline{N}\ra \Gamma$
is the canonical projection, and $j$ is a homotopy inverse of the
equivalence
$B\overline{N}^{\wedge}_p\stackrel{Bi_{\overline{N},\overline{G}}}{\simeq
} B\overline{G}^{\wedge}_p$. As $B\Gamma ^{\wedge}_p$ is already
$p$-complete, the universal property of $p$-completion show that
is enough to check that $g\circ f\circ (-)^{\wedge}_p$ is
homotopic to $(-)^{\wedge}_p.$ Now, applying Corollary
\ref{corozabrodsky} to the fibre sequence $B\overline{A}\ra
B\overline{N} \xrightarrow{B\rho}B\Gamma$ (with
$X=B\Gamma^{\wedge}_p$) we obtain that $g\circ f\circ
(-)^{\wedge}_p\simeq (-)^{\wedge}_p$ if and only if $g\circ f\circ
(-)^{\wedge}_p\circ B\rho \simeq (-)^{\wedge}_p\circ B\rho$. Again
because $B\Gamma^{\wedge}_p$ is $p$-complete, the latter is
equivalent to prove that the map $B\rho^{\wedge}_p$ is homotopic
to $g\circ f\circ B\rho^{\wedge}_p$. Again by completeness and
Corollary \ref{corozabrodsky} applied to the fibration
$B\O_A(G)\ra BG \xrightarrow{B\pi} B\overline{G}$, it is seen that
we should check that $B\rho^{\wedge}_p\circ j\circ B\pi\simeq
B\rho^{\wedge}_p\circ j\circ B\pi\circ g\circ f$.

Now it is clear by construction of $\alpha$ that
$B\rho^{\wedge}_p\circ j\circ B\pi^{\wedge}_p$ is homotopic to
$\alpha$. Therefore, we need to see that $g\circ
f\circ\alpha\simeq\alpha.$ By diagram \ref{diag6} the latter is
homotopic to $g\circ\eta\circ h$, and has the same homotopy class
as $\alpha$ by diagram~\ref{diag1}. So we are done.
\end{proof}

When we combine the previous statement with the last statement of
Proposition~\ref{cofibre} we obtain a complete description of
$\mathbf{CW}_{B\Z /p}BG$ for every $p$ and every finite group $G$.

\begin{theorem} Let $G$ be a finite group generated by its elements of order
$p$, let $S \in Syl_p(G)$, and let $A = \omegabar S$ be the
minimal strongly closed subgroup of $S$ containing $\Omega_1(S)$.
Let overbars denote passage from $G$ to the
quotient group $G / \O_A(G)$.
Then the $B\Z/p$-cellularization of $BG$ has one the following
shapes:

\begin{itemize}

\item[(1)]
If $G = S$ is a $p$-group then $BG$ is $B\Z /p$-cellular.

\item[(2)]
If $G$ is not a $p$-group and $A = S$ then $\mathbf{CW}_{B\Z/p}BG$
is the homotopy fiber of the natural map $BG \ra \prod_{q\neq
p}BG^{\wedge}_q$.

\item[(3)]
If $G$ is not a $p$-group and $A \ne S$ then
$\mathbf{CW}_{B\Z/p}BG$ is the homotopy fiber of the map $BG \ra
B(N_{\overline G}(\overline A)/\, \overline A)^{\wedge}_p \times
\prod_{q\neq p}BG^{\wedge}_q$.
\end{itemize}
\label{recopcel}
\end{theorem}

\begin{remark}

It is worth recalling here that if $G$ is not generated by order
$p$ elements, its $B\Z /p$-cellularization can be computed by
using the aforementioned homotopy equivalence
$\mathbf{CW}_{B\Z/p}B\Omega_1(G)\simeq\mathbf{CW}_{B\Z/p}BG$
induced by inclusion, and then applying the previous theorem to
$\Omega_1(G)$.

\end{remark}

Theorems 1.1 and 1.2 and Corollary 1.5 in \cite{Flores08}
determine $N_G(\overline A)/\overline A$, whose structure is very
rigid and depends on a restricted set of well-known simple groups.
It is very likely that an analogous classification can be obtained
exactly in the same way for $\mathbf{CW}_{B\Z /p^r}BG$, $r>1$, but
we have restricted ourselves to the case $r=1$ for the sake of
simplicity (cf. also \cite[Theorem 3.6]{CCS2}).

\medskip

In the cases where $\O_A(G)=1$ and $A\nor G$
--- which are implicit in the computations --- the $B\Z
/p$-cellularization of $BG^{\wedge}_p$ is the homotopy fiber of
the natural map $BG^{\wedge}_p\ra B(G/A)^{\wedge}_p$. It is then
tempting to identify $\mathbf{CW}_{B\Z /p}BG$ with $BA$. But this
would mean, in particular, that
$\textrm{map}_*(BZ/p,BG^{\wedge}_p)$ would be discrete. However,
an analysis of the fibration sequence
$$\textrm{map}_*(B\Z /p,BG^{\wedge}_p)\ra\textrm{map}(B\Z /p,BG^{\wedge}_p)\ra
BG^{\wedge}_p$$ together with the description of its total space
--- which is given, for example, in \cite[Appendix]{Broto02}
--- shows that
$\textrm{map}_*(B\Z /p,BG^{\wedge}_p)$ is non-discrete in
general, and then usually $\mathbf{CW}_{B\Z /p}BG$ is not an
aspherical space.

It is conceivable that our results can also have interesting
consequences from the point of view of homotopical representations
of groups. In \cite[Section 6]{Flores07} the results on
cellularization gave rise to specific examples of nontrivial maps
$BG\ra BU(n)^{\wedge}_p$ that enjoyed two particular properties:
they did not come from group homomorphisms $G\ra U(n)$, and they
were trivial when precomposing with any map $B\Z /p\ra BG$. While
there are a number of examples in the literature with the first
feature (see for example \cite{Benson95} or \cite{Mislin89}), no
representations were known at this point for which the second
property holds. The classification results of this paper give hope
of finding a systematic and complete treatment of all these kinds
of representations. We plan to undertake this task in a separate
paper.

In the next section we show the applicability of our results
by computing the $B\Z /p$-cellularization of various specific
families of classifying spaces. We have chosen the simple groups
(as they have shown their cornerstone role in the computation of
$\mathbf{CW}_{B\Z /p}BG$), certain split extensions that signaled
there was something beyond the results of
\cite{Flores07}, and certain nonsplit extensions of $G_2(q)$ that
illuminate the roles of the normalizers of $A$ and $S$
in the $B\Z /p$-cellular context.

\section{Examples}
\label{examples}

In Theorem~\ref{recopcel} a description of the
$B\Z/p$-cellularization of $BG$ for every group $G$ and every
prime $p$ is given, so in this section we describe some families
of concrete examples for which $\mathbf{CW}_{B\Z /p}B G$ is
interesting. We begin by providing more explanation and some examples
of the group-theoretic concepts.

As observed earlier, if $N$ is any normal subgroup of $G$, the
Sylow $p$-subgroup $S \cap N$ of $N$ is strongly closed in $G$.
More generally, if $A_0$ is a strongly closed $p$-subgroup of a group $G_0$
and $\pi : G \rightarrow G_0$ is any surjective homomorphism
with kernel $N$, then any Sylow $p$-subgroup, $A$, of the complete
preimage $\pi^{-1}(A_0)$ is strongly closed in $G$ with
$N \le \O_{A}(G)$.  Since such extensions of a given $G_0$
are essentially arbitrary,
it is natural to seek the ``minimal'' groups possessing strongly
closed subgroups, $A$, namely those with $\O_A(G) = 1$.
This is achieved by passing to the quotient $\overline G = G/\O_A(G)$
wherein $\overline A$ is still strongly closed (but possibly trivial)
and does not contain the Sylow $p$-subgroup of any nontrivial normal subgroup.
The main result of \cite{Flores08} is a characterization of all these
``minimal configurations.'' For the reader's convenience we now state a
slightly weaker but less recondite version of that classification for all primes $p$.

\begin{theorem}\label{strongly-closed}
Let $A$ be any nontrivial strongly closed $p$-subgroup of $G$.
Assume $\O_A(G) = 1$ and $G$ is generated by conjugates of $A$.
Then $G$ has a normal subgroup $L$ such that
$$
L \iso L_1 \times L_2 \times \cdots \times L_r
$$
where each $L_i$ is a simple group, $A_i = A \cap L_i$ is strongly
closed but not Sylow in $L_i$ and
$A = (A_1 \times \cdots A_r)A_F$ where the subgroup $A_F$ (possibly trivial)
normalizes each $L_i$ and induces field automorphisms on $L_i$.
Each $L_i$ is one of the following types:
(i) a Chevalley group of arbitrary
$BN$-rank in characteristic $\ne p$ with
abelian but not elementary abelian Sylow $p$-subgroups,
(ii) a Chevalley group of $BN$-rank~1 in characteristic $= p$,
(iii) $G_2(q)$ with $(q,3) = 1$ and $p = 3$, or
(iv) one of seven sporadic groups (with $p$ specified in each case).
\end{theorem}

Here $L$ is the socle of $G$ (also $L = F^*(G)$) and the exact structure of $G$
is described in \cite{Flores08}.
Moreover, the full classification describes all possibilities for $A_i$,
its normalizer, and how the conditions in family (i)
may easily be determined.  The ``field automorphism subgroup''
$A_F$ of $A$ can only be nontrivial when there are factors from family (i).

In the special case when $A = \omegabar S$, the condition that
$G$ be generated by its elements of order $p$ is equivalent to $G$
being generated by conjugates of $A$.  The classification shows that
when $A = \omegabar S \ne S$ and
$\O_A(G) = 1$, then $A = \Omega_1(S)$ and $A$ is elementary abelian
(although $S$ itself need not be abelian).
We shall see examples in Subsection~\ref{split-extensions}
where $\omegabar S$ is not abelian when $\O_A(G) \ne 1$.

In our context, two important (slightly simplified renditions of) consequences
of this theorem are:

\begin{corollary}
\label{control-fusion}
Under the hypotheses of Theorem~\ref{strongly-closed}

\begin{itemize}
\item[(1)]
$N_G(A)$ controls strong fusion in $S$.
\item[(2)]
If $A_F = 1$ then $\Gamma = N_G(A)/A = N_{L_1}(A_1)/A_1 \times \cdots \times N_{L_r}(A_r)/A_r$.
\end{itemize}
\end{corollary}

When $A_F \ne 1$ the general classification also describes $\Gamma$ precisely---it may be viewed as a
subgroup of the direct product in (2); and in all cases $\Gamma$ is $p$-perfect.

One example where $L = L_1$ belongs to family (i) but
$A_F \ne 1$ is when $G = PSL_{11}(q)\gp{f}$ with $p=5$, $q = 3^5$, and
$f$ is a field automorphism of order 5 (this example is explicated
in greater detail in \cite{Flores08}).
Here the simple group
$L = PSL_{11}(q)$ has an abelian Sylow 5-subgroup of type (25,25).
If $f \in S \in Syl_5(G)$, then
$$
A = \omegabar S = \Omega_1(S) = \gp{f,\Omega_1(S \cap L)}
\iso Z_5 \times Z_5 \times Z_5
$$
and $A$ is evidently strongly closed in $S$ with respect to $G$.
Furthermore $A^* = \Omega_1(S \cap L)$ is a strongly closed
subgroup of $L$. Both $G$ and $L$ are generated by their elements of order 5.
One may also calculate that
$$
\Gamma = N_G(A)/A \iso
(H \times H)\gp t \times GL_3(3)
$$
where $H$ is a split extension $Z_{80}\cdot Z_4$ which has no
subgroup of index 5, and $t$ has order 2, commutes with $f$ and
interchanges the two copies of $H$.
Thus $\Gamma$ has no subgroup of index 5, even when $G$ itself has a normal subgroup of index 5.

In Subsections~\ref{split-extensions} and \ref{exotic} we explore control of
fusion in extensions of the simple groups possessing strongly closed subgroups:
in particular, we examine when $N_G(S)$, which is a subgroup of $N_G(A)$,
also controls fusion in $S$.

\subsection{Simple groups}
\label{simple-groups}

In the previous discussion we sketched the list of all
simple groups $G$ possessing a strongly closed $p$-subgroup $A$ that is not
a Sylow $p$-subgroup. By simplicity, $G$ is generated by
its elements of order $p$, and $\O_A(G) = 1$.
In \cite[5.6-5.8]{Flores07} the $B\Z /2$-cellularization of the
classifying spaces of all the simple groups was computed.
(For $p=2$ only the groups $U_3(2^n)$ and $Sz(2^n)$ can occur;
these belong to family (ii).)
In this section we invoke the classification for $p$ odd,
Theorem~\ref{strongly-closed},
to see how for
every simple group $G$ the cellularization $\mathbf{CW}_{B\Z /p}B
G$ is included in the different cases of Theorem~\ref{recopcel}.
By this classification, when $G$ is simple and $\omegabar S \ne S$
we always have $N_G(\omegabar S) = N_G(S)$ (see Corollary~2.8 in \cite{Flores08}).
We thus obtain the following characterization:

\begin{proposition}

Let $G$ be a simple group,
let $p$ a prime and let $S$ be a Sylow $p$-subgroup of $G$.
Then $\mathbf{CW}_{B\Z /p}BG$ has one of the following two structures:

\begin{itemize}

\item[(1)]
If $\omegabar S = S$, then $\mathbf{CW}_{B\Z /p}BG$ is the
homotopy fiber of the natural map $BG\ra \prod_{q\neq p}BG^{\wedge}_q.$

\item[(2)]
If $\omegabar S \ne S$, then we have a fibration sequence
$$\mathbf{CW}_{B\Z /p}BG\ra BG\ra
 B(N_{G}(S)/\omegabar S)^{\wedge}_p\times\prod_{q\neq
p}BG^{\wedge}_q.$$

\end{itemize}

\end{proposition}

Note that the inclusion $N_G(S)\hookrightarrow N_G(\omegabar S)$
induces a homotopy equivalence $BN_G(S)^{\wedge}_p\simeq
BN_G(\omegabar S)^{\wedge}_p$ when $N_G(S)$ (and then
$N_G(\omegabar S)$) controls $p$-fusion in $S$. This happens ($p$
odd) for every simple group such that $\omegabar S \ne S$, and in
particular a comparison of the fibration sequences of (Theorem
\ref{Chacholski}) for $A=B\Z /p$ and $X=BN_G(S)^{\wedge}_p$ and
$X=BN_G(\omegabar S)^{\wedge}_p$ (which are homotopy equivalent)
respectively, gives that the induced map $B(N_G(S)/\omegabar
S)^{\wedge}_p\simeq B(N_G(\omegabar S)/\omegabar S)^{\wedge}_p$ is
also a homotopy equivalence.

The general structure of the normalizers that appear in the previous
characterization is described in Sections~2 and 4.1 of \cite{Flores08}.
For a specific example, take $p$ odd and let $q$ be any power of an odd prime
such that $p^2 \divides q-1$,
so that a Sylow $p$-subgroup $S$ of $G = PSL_2(q)$ is cyclic of order $\ge p^2$
(for example, $p=3$ and $q = 19$).
One easily computes that for $A = \omegabar S = \Omega_1(S) < S$
we have $N_G(A) = N_G(S)$ is dihedral of order
$q-1$ and so $N_G(A)/A$ is dihedral of order $(q-1)/p$.

\medskip

\subsection{Split extensions}
\label{split-extensions}

We turn now to the case of non-simple groups. Here we give some
explicit examples of $B\Z /p$-cellularization of split extensions
that are beyond the scope of \cite{Flores07}, and show the
usefulness of Theorems \ref{thetheorem} and \ref{recopcel}.

In \cite{Flores07} the $B\Z /p$-cellularization of $BG$ is
described when $G$ is generated by elements of order $p$,
$\omegabar S$ is a proper subgroup of $S$, and the normalizer of
$S$ controls strong fusion in $S$. No example was given there of a
group for which the first two conditions hold but not the third.
George Glauberman suggested an example of a group of the latter
type: a wreath product $(\Z /2) \wr Sz(2^n)$. In this section we
generalize this example, showing that many split extensions for
which these conditions hold can be constructed. The computation of
the cellularization of their classifying space is then easy from
Corollaries 1.4 and 1.5 in \cite{Flores08}.

Let $R$ be a simple group, and assume that there exists a Sylow
$p$-subgroup $T$ of $R$ for which $\omegabar T \ne T$ (for instance
any simple group in family (i) of Theorem~\ref{strongly-closed}).
Let $E$ be any elementary abelian $p$-group on which $R$ acts faithfully.
Then $S=ET$ is a Sylow $p$-subgroup of the semi-direct product
$G=E \sdp R$, and if we denote $A = E\omegabar T$, it is clear
that $E = \O_A(G)$. As the extension is split, the canonical
projection $G \ra R$ sends $\omegabar S$ to $\omegabar T$.
Note that because $\omegabar T$ acts faithfully on $E$,
$\omegabar S$ is non-abelian.
In such a group $G$ {\it neither $N_G(A)$ nor $N_G(S)$ controls
fusion in $S$} (see Section~4.2 of \cite{Flores08} for details).
However, in the quotient group $\overline G \iso R$
the normalizer $N_{\overline G}(\overline A)$ does control fusion
in $\overline S$.  This further illuminates why the classification
is needed to pass to the ``correct'' quotient for our fusion arguments.

For a specific example, take $R = PSL_2(q)$
satisfying the conditions at the end of the previous subsection
(where it was called $G$) and take $E$ to be the $\FF_p R$-module
affording the regular representation of $R$ over the field $\FF_p$.

\medskip

According to Theorem~\ref{thetheorem}, the $B\Z
/p$-cellularization of $BG^{\wedge}_p$ has the homotopy type of
the fiber of the composite
$$
BG^{\wedge}_p\longrightarrow BR^{\wedge}_p \longrightarrow
B(N_{R}(\omegabar T)/\omegabar T)^{\wedge}_p.
$$

The structure of $G$ was studied in the \cite[4.1]{Flores08}, and
moreover, Section~4.2 there also describes the shape of the
normalizer in a concrete example where $R$ is a special linear
group and $T$ is a subgroup of diagonal matrices.


\subsection{Exotic extensions of $\pmb{G_2(q)}$}
\label{exotic}

We conclude with somewhat exotic examples which show that in
Theorem \ref{recopcel}, the normalizer of the subgroup
$\overline{A}$ in $\overline{G}$ cannot be replaced by the normalizer of the
corresponding Sylow subgroup. These groups are constructed as
``half-split" extensions of the simple group $G_2(q)$, for $q$
coprime to 3, as established in \cite{Flores08}.
(The simple groups $G_2(q)$ are in family (iii) of Theorem~\ref{strongly-closed}.)
For convenience we restate Theorem~4.4 of \cite{Flores08} here
(whose proof is constructive):

\begin{proposition}
Let $p$ be a prime dividing the order of the finite group $R$
and let $X$ be a subgroup of order $p$ in $R$. Then there is an
$\FF_p R$-module $E$ and an extension
$$
1 \longrightarrow E \longrightarrow  G \longrightarrow  R \longrightarrow  1
$$
of $R$ by $E$ such that the
extension of $X$ by $E$ does not split, but the extension of $Z$ by
$E$ splits for every subgroup $Z$ of order $p$ in $R$ that is not conjugate
in $R$ to $X$.
In particular, for nonidentity elements $x \in X$ and $z \in Z$ every element
in the coset $xE$ has order $p^2$ whereas $zE$ contains elements of order $p$
in $G$.

\end{proposition}

We study the case $R = G_2(q)$ with $p = 3$ and $(q,3) = 1$, and
we let $T \in Syl_3(R)$. The description in
\cite[Proposition~2.7(3)]{Flores08} of the normalizer of $T$ in $R$ implies that
$BN_R(T)^{\wedge}_3$ is not homotopy equivalent to
$BG_2(q)^{\wedge}_3$. Now apply the previous Proposition with
$Z = Z(T)$ and $X = \gp x$ for any
$x \in T - Z$ of order 3 (all such subgroups $X$ are conjugate in $G$,
but none are conjugate to $Z$).  Then if $S$ is the Sylow 3-subgroup
of $G$ containing $E$ such that $S/E = T$, then the ``half-split''
construction forces $E \le \omegabar S$ and $\omegabar S /E = ZE/E$
(unlike a split extension of $R$ by $E$ where $\omegabar S = S$).
Again the considerations in \cite[Section~4.3]{Flores08} and
Theorem~\ref{thetheorem} establish that the
$B\Z /3$-cellularization of $(BG)^{\wedge}_3$ has the homotopy
type the homotopy fibre of the map $(BG)^{\wedge}_3\ra
BPSL^*_3(q)^{\wedge}_3$, but not of the map $(BG)^{\wedge}_3\ra
B(N_R(T)/T)^{\wedge}_3$ when $9 \divides q^2 - 1$. Here
$N_R(Z) \iso PSL^*_3(q)$, where the latter group denotes the
projective version of the
group $SL_3(q)$ or $SU_3(q)$ together with the outer (graph)
automorphism of order 2 inverting its center.

From this computation one can deduce, then, that the object that
determines the $\Sigma B\Z /p$-nullification of the cofibre of the
map of Theorem \ref{Chacholski} in the case of finite groups is
the normalizer of $\overline{A}$, and not the Sylow normalizer, as
might be inferred from the particular cases studied in
\cite{Flores07}. This example also highlights the importance of
having a classification of \emph{all} groups possessing a
nontrivial strongly closed $p$-subgroup that is not Sylow --- not
just the simple groups having such a subgroup that contains
$\Omega_1(S)$
--- since the subgroup $\omegabar S$ does not pass in a
transparent fashion to quotients.

\medskip

In conclusion, some interesting open questions remain. We have
characterized with precision $\mathbf{CW}_{B\Z /p}BG$ for every
finite $G$, and in the course of the proof we have also described
$\mathbf{CW}_{B\Z /p}BG^{\wedge}_p$ when $G$ is generated by
order $p$ elements. However we do not address the issue of what
happens in general with the cellularization of $BG^{\wedge}_p$ if
we remove the generation hypothesis. There are some cases that can
be deduced from the previous developments --- for example if $G$
is not equal to $\Omega_1(G)$ but is mod $p$ equivalent to a group
that is so --- but it would be nice to have a general statement.

The extensions of our techniques and results to more general
$p$-local spaces with a notion of $p$-fusion seem to be the
natural next step of our study; in particular, classifying spaces
of $p$-local finite groups and some families of non-finite groups
offer enticing possibilities.


\newpage





\bigskip

\bigskip\noindent
Ram\'on J.~Flores\\
Departamento de Estad\'\i stica, Universidad Carlos III de
Madrid,\\
Avda. de la Universidad Carlos III, 22\\
E -- 28903 Colmenarejo (Madrid) --- Spain \\
e-mail: {\tt rflores@est-econ.uc3m.es}

\medskip

\noindent
Richard M.~Foote\\
Department of Mathematics and Statistics, University of Vermont, \\
16 Colchester Avenue,
Burlington, Vermont 05405 --- U.S.A. \\
e-mail: {\tt foote@math.uvm.edu}




\begin{thebibliography}{draft7}

\bibitem[AM94]{MR1317096}
A.~Adem and R.~J.~Milgram.
\newblock {\em Cohomology of finite groups}, volume 309 of {\em Grundlehren der
  Mathematischen Wissenschaften [Fundamental Principles of Mathematical
  Sciences]}.
\newblock Springer-Verlag, Berlin, 1994.

\bibitem[BK72]{Bousfield72}
A.~K.~Bousfield and D.~M.~Kan.
\newblock {\em Homotopy limits, completions and localizations}.
\newblock Springer-Verlag, Berlin, 1972.
\newblock Lecture Notes in Mathematics, Vol. 304.

\bibitem[BK02]{Broto02}
C.~Broto and N.~Kitchloo.
\newblock Classifying spaces of Kac-Moody groups.
\newblock {\em Math. Z.}, 240(2):621--649, 2002.

\bibitem[BKI08]{BKI08}
D.~Benson, H.~Krause and S.~Iyengar.
\newblock Local cohomology and support for triangulated
categories.
\newblock To appear in {\em Ann. Sci {\'E}cole Norm. Sup}.

\bibitem[BLO03]{Broto032}
C.~Broto, R.~Levi, and B.~Oliver.
\newblock The homotopy theory of fusion systems.
\newblock {\em J. Amer. Math. Soc.}, 16(4):779--856 (electronic), 2003.

\bibitem[BLO07]{Broto07}
C.~Broto, R.~Levi, and B.~Oliver.
\newblock Discrete models for
the $p$--local homotopy theory of compact Lie groups and
$p$--compact groups.
\newblock {\em Geometry and Topology}, 11:315--427, 2007.

\bibitem[Bou75]{Bousfield75}
A.K. Bousfield.
\newblock The localization of spaces with respect to homology.
\newblock {\em Topology}, 14:133--150, 1975.

\bibitem[Bou94]{Bousfield94}
A.K.~Bousfield.
\newblock Localization and periodicity in unstable homotopy theory.
\newblock {\em J. Amer. Math. Soc.}, 7(4): 831--873, 1994.

\bibitem[BW95]{Benson95}
D.J.~Benson and C.W.~Wilkerson.
\newblock Finite simple groups and {D}ickson invariants.
\newblock In {\em Homotopy theory and its applications (Cocoyoc, 1993)}, volume
  188 of {\em Contemp. Math.}, pages 39--50. Amer. Math. Soc., Providence, RI,
  1995.

\bibitem[CCS07]{CCS2}
N.~Castellana, J.A.~Crespo, and J.~Scherer.
\newblock {P}ostnikov pieces and ${B}{\mathbb {Z}}/p$-homotopy theory.
\newblock {\em Trans. Amer. Math. Soc.}, 359(4):1791--1816, 2007.

\bibitem[CDI06]{CDI06}
W.~Chach{\'o}lski, W.G.~Dwyer and M.~Intermont.
\newblock The {$A$}-complexity of a space.
\newblock {\em J. London Math. Soc. (2)}, 65(1):204-222, 2006.

\bibitem[Cha96]{Chacholski96}
W.~Chach{\'o}lski.
\newblock On the functors {$CW\sb A$} and {$P\sb A$}.
\newblock {\em Duke Math. J.}, 84(3):599--631, 1996.

\bibitem[DGI01]{Dwyer06}
W.G.~Dwyer, J.P.C.~Greenlees and S.~Iyengar.
\newblock Duality in algebra and topology.
\newblock {\em Adv. Math. 200}, 200(2):357--402, 2006.

\bibitem[DGS07]{DGS07}
\newblock E. Dror-Farjoun, R.~G{\"o}bel and Y.~Segev.
\newblock Cellular covers of groups.
\newblock {\em J. Pure Appl. Algebra}, 208(1):61--76, 2007.

\bibitem[DGS08]{DGS08}
\newblock E. Dror-Farjoun, R.~G{\"o}bel, Y.~Segev and S.~Shelah.
\newblock On kernels of cellular covers.
\newblock {\em Groups Geom. Dyn. 1}, 4:409--419, 2007.

\bibitem[Dwy96]{Dwyer96}
W.G.~Dwyer.
\newblock The centralizer decomposition of {$BG$}.
\newblock In {\em Algebraic topology: new trends in localization and
  periodicity (Sant Feliu de Gu\'\i xols, 1994)}, volume 136 of {\em Progr.
  Math.}, pages 167--184. Birkh\"auser, Basel, 1996.

\bibitem[Far95]{Farjoun95}
E.~Dror-Farjoun.
\newblock {\em Cellular spaces, null spaces and homotopy localization}, volume
  1622 of {\em Lecture Notes in Mathematics}.
\newblock Springer-Verlag, Berlin, 1996.

\bibitem[FGS07]{Farjoun07}
E.~Dror-Farjoun, R.~G\"obel and Y.~Segev.
\newblock Cellular covers of spaces.
\newblock To appear in {\em J. Pure Appl. Algebra}.

\bibitem[Flo07]{Ramon}
R.J.~Flores.
\newblock Nullification and cellularization of classifying spaces of finite
  groups.
\newblock {\em Trans. Amer. Math. Soc.}, 359, 1791--1816, 2007.

\bibitem[FF08]{Flores08}
R.J.~Flores and R.~Foote.
\newblock Strongly closed subgroups of finite groups.
\newblock {\em Advances in Math.}, 222, 453--484, 2009.


\bibitem[Foo97]{Foote97}
R.~Foote.
\newblock A characterization of finite groups containing a strongly closed
  {$2$}-subgroup.
\newblock {\em Comm. Algebra}, 25(2):593--606, 1997.

\bibitem[FS07]{Flores07}
R.J.~Flores and J.~Scherer.
\newblock Cellularization of classifying spaces and fusion properties of finite groups.
\newblock {\em J. Lond. Math. Soc.(2)} 76(1):41--56, 2007.

\bibitem[Gol74]{Gol74}
D.~Goldschmidt.
\newblock 2-Fusion in finite groups.
\newblock {\em Ann. Math.}, 99(1974), 70--117.

\bibitem[Gol75]{Gol75}
D.~Goldschmidt.
\newblock Strongly closed 2-subgroups of finite groups.
\newblock {\em Ann. Math.}, 102(1975), 475--489.

\bibitem[Kie08]{Kiessling08}
J.~Kiessling.
\newblock Classification of certain cellular classes of chain
complexes.
\newblock Preprint, available at arXiv:0801.3904.

\bibitem[MP98]{MP98}
J.~Martino and S.~Priddy.
\newblock On the cohomology and homotopy of Swan groups.
\newblock {\em Math. Z.}, 225(2):277--288, 1997.
  RI, 1998.

\bibitem[Mil84]{Miller}
H.~Miller.
\newblock The {S}ullivan conjecture on maps from classifying spaces.
\newblock {\em Ann. of Math. (2)}, 120(1):39--87, 1984.

\bibitem[MT89]{Mislin89}
G.~Mislin and C.~Thomas.
\newblock On the homotopy set {$[B\pi,BG]$} with {$\pi$} finite and {$G$} a
  compact connected {L}ie group.
\newblock {\em Quart. J. Math. Oxford Ser. (2)}, 40(157):65--78, 1989.

\bibitem[RS01]{Rodriguez01}
J.~L.~Rodr{\'{\i}}guez and J.~Scherer.
\newblock Cellular approximations using {M}oore spaces.
\newblock In {\em Cohomological methods in homotopy theory (Bellaterra, 1998)},
  volume 196 of {\em Progr. Math.}, pages 357--374. Birkh\"auser, Basel, 2001.

\bibitem[RS08]{RS08}
\newblock J.~L.~Rodr{\'{\i}}guez and J.~Scherer.
\newblock A connection between cellularization for groups and spaces via
two-complexes.
\newblock {\em J. Pure Appl. Algebra}, 212:1664-1673, 2008.



\end{thebibliography}
\end{document}